\documentstyle[12pt]{article}

\begin{document}
\begin{center}
{\large\bf
A combinatorial characterization of second category subsets
of $X^\omega$}
\end{center}
\vskip 1.0truecm
\centerline{{\bf Apoloniusz Tyszka}}
\vskip 1.0truecm
\def\thefootnote{}
\begin{abstract}
\noindent
Let a finite non-empty $X$ is equipped with discrete topology.
We prove that $S \subseteq X^\omega$ is of second category
if and only if for each
$f:\omega\rightarrow\bigcup_{n\in\omega}X^n$ there exists
a sequence $\{a_n\}_{n\in\omega}$ belonging to $S$ such that for
infinitely many $i\in\omega$ the infinite sequence
$\{a_{i+n}\}_{n\in\omega}$ extends the finite sequence $f(i)$.
\end{abstract}
\vskip 1.0truecm
\par
Theorem~1 yields information about sets $S\subseteq X^\omega$
with the following property $(\Box)$:
\begin{description}
\item[$(\Box)$]
for each infinite $J\subseteq \omega$ and each
$f:J\rightarrow\bigcup_{n\in\omega}X^n$ there
exists a sequence $\{a_n\}_{n\in\omega}$ belonging to $S$ such
that for infinitely many $i\in J$ the infinite sequence
$\{a_{i+n}\}_{n\in\omega}$ extends the finite sequence $f(i)$.
\end{description}

\par
{\bf Theorem~1.}
Assume that a non-empty $X$ is equipped with discrete topology.
We claim that if $S \subseteq X^\omega$
is of second category then $S$ has the property $(\Box)$.
\par
{\it Proof.} Let us fix
$f:J\rightarrow\bigcup_{n\in\omega}X^n$. Let $S_k (f)$
$(k\in\omega)$
denote the set of all sequences $\{a_n\}_{n\in\omega}$
belonging to $X^\omega$ with the property that there exists
$i\in J$ such that $i>k$ and the infinite sequence
$\{a_{i+n}\}_{n\in\omega}$ extends the finite sequence $f(i)$.
Sets $S_k (f)$ $(k\in\omega)$ are open and dense. In virtue of the
Baire category theorem
$\bigcap_{k\in\omega} S_k (f)\cap S$ is non-empty i.e.
there exists a sequence $\{a_n\}_{n\in\omega}$ belonging to $S$
such that
for infinitely many $i\in J$ the infinite sequence
$\{a_{i+n}\}_{n\in\omega}$ extends the finite sequence $f(i)$.
This completes the proof.
\par
\vskip 0.3 truecm
The proof of the following Observation is left as an exercise for
the reader.
\vskip 0.3truecm
{\bf Observation}. If $S \subseteq \{0,1\}^\omega$ has the property
$(\Box)$ then for every open set $U \subseteq (0,\varepsilon)$ with
$0 \in \overline{U}$ there exists a $g \in S$ such that the sequence
\begin{displaymath}
\{\sum_{k=n}^\infty \frac {g(k)}{2^k}\}_{n \in \omega}
\end{displaymath}
has an infinite number of terms belonging to $U$.
\vskip 0.3truecm
{\bf Corollary.} Assume that
$f:(0,\varepsilon) \rightarrow {\bf R}$ is
continuous and for each zero-one sequence $\{a_n\}_{n \in \omega}$
with an infinite number of ones the limit
\begin{displaymath}
\lim_{n \rightarrow \infty}
f(\sum_{k=n}^\infty \frac {a_k}{2^k})
\end{displaymath}
exists and equals $0$.
Then (cf. Proposition 1 in [5]) $\lim_{x \rightarrow 0^+} f(x)=0$.
\vskip 0.5truecm
Theorem~2 yields information about sets $S\subseteq X^\omega$
with the following property $(\ast)$:
\begin{description}
\item[$(\ast)$]
for each $f:\omega\rightarrow\bigcup_{n\in\omega}X^n$ there
exists a sequence $\{a_n\}_{n\in\omega}$ belonging to $S$ such that
for infinitely many $i\in\omega$ the infinite sequence
$\{a_{i+n}\}_{n\in\omega}$ extends the finite sequence $f(i)$.
\end{description}

\par
{\bf Theorem~2.}
Assume that a finite non-empty $X$ is equipped with discrete
topology. We claim that if $S \subseteq X^\omega$ is of
first category then $S$ does not have the property $(\ast)$.
\par
{\it Proof.}
Assume that
$S \subseteq\bigcup_{i \in \omega} Y_{i}$ , where sets
$Y_{0} \subseteq Y_{1} \subseteq Y_{2} \subseteq... X^\omega$
are closed and nowhere dense.
Let $\Pi_{i}:X^\omega \rightarrow X^\omega$
($i\in\omega$)
maps the sequence $\{a_{n}\}_{n \in \omega}\in X^\omega$ to the
sequence $\{a_{i+n}\}_{n \in \omega}$. Obviously
$\Pi_{i}=
\underbrace{\Pi_{1}\circ...\circ\Pi_{1}}_{i \; {\rm times}}$
($i\in\omega\setminus\{0\}$). Since $X$ is finite,
$\Pi_{1}$ images of closed nowhere dense sets are closed and
nowhere dense.
Therefore each set $\Pi_{i}(Y_{i})$ ($i\in\omega$) is closed
and nowhere dense.
Hence for each $i\in\omega$ there exists a sequence
$b_i(0),\:b_i(1),\:.\:.\:.\:,\:b_i(l(i))$ of elements of $X$
such that the set

\vskip 0.2truecm
\centerline{$\{b_i(0)\}\times\{b_i(1)\}\times...\times\{b_i(l(i))\}
\times X\times X\times X\times...$}

\vskip 0.2truecm
\noindent
is disjoint from $\Pi_{i}(Y_{i})$.
This gives:
\begin{description}
\item[$(\ast\ast)$]
if $i \in \omega$ then each
sequence $\{a_n\}_{n\in\omega} \in X^\omega$ which satisfies
$a_{i}=b_{i}(0),a_{i+1}=b_{i}(1),...,a_{i+l(i)}=b_{i}(l(i))$
does not belong to $Y_{i}$.
\end{description}

\noindent
Let $f(i)$ ($i\in\omega$) denote the sequence
$b_{i}(0),\;b_{i}(1),\;.\;.\;.\;,\;b_{i}(l(i))$; formally
$f:\omega\rightarrow\bigcup_{n\in\omega}X^n$.
Let $\{a_n\}_{n\in\omega}\in S$ and
$I:=\{i\in\omega: \{a_{i+n}\}_{n\in\omega}\ {\rm extends}\ f(i) \}$.
Obviously $I \subseteq \omega$, it suffices to prove that $I$
is finite. From $(\ast\ast)$ we conclude that for each
$i \in I \;\; \{a_n\}_{n\in\omega} \not\in Y_{i}$.
Suppose, on the contrary, that $I$ is infinite.
Thus $\{a_n\}_{n\in\omega} \not\in \bigcup_{i \in I}
Y_{i}=\bigcup_{i\in\omega} Y_{i}$.
Since $S \subseteq \bigcup_{i\in \omega} Y_{i}$ we conclude that
$\{a_n\}_{n\in\omega}\not\in S$, which contradicts our assumption.
We have proved that $S$ does not have the property $(\ast)$.
\vskip 0.3truecm
{\bf Remark~1} ([3]). Tomek Bartoszy\'{n}ski constructed a closed
nowhere dense set $S \subseteq \omega^\omega$ with the property
$(\ast)$.

{\bf Remark~2} (inspired by [3]). Let $X$ is infinite,
$\psi: X \rightarrow \omega$ and $\psi(X)$ is infinite. Let
$S \subseteq X^\omega$ denote the set of all sequences of the form
\begin{displaymath}
a,\underbrace{.........................}_{\psi(a) \ {\rm elements}
 \ {\rm of} \ X},a,
b,\underbrace{.........................}_{\psi(b) \ {\rm elements}
 \ {\rm of} \ X},b,
c,\underbrace{.........................}_{\psi(c) \ {\rm elements}
 \ {\rm of} \ X},c,
...
\end{displaymath}
where
$a,b,c,... \in X$. It is easy to check that $S$ is closed, nowhere
dense and has the property $(\Box)$.

\vskip 0.3truecm

Let $\forall^\infty$ abbreviate "for all except finitely many".

{\bf Note.} If $S \subseteq \omega^\omega$ satisfies
$\exists g \in \omega^\omega\ \forall f \in S \ \forall^\infty k \ \
g(k) \neq f(k)$ then $S$ does not have the property $(\ast)$.
\vskip 0.3truecm
Let ${\cal M}$ denote the ideal of first category subsets of
{\bf R} and ${\rm non}({\cal M}):
={\rm min}\{{\rm card}\ S: S \subseteq {\bf R},\; S \not\in
{\cal M}\}$.
It is known (see [1], [2] and also [4]) that:
\vskip 0.2truecm
\noindent
\centerline{
non$(\cal M)$ = 
min$\{{\rm card}\ S: S\subseteq\omega^\omega\ {\rm and} \
\neg \ \exists g \in \omega^\omega \
\forall f \in S \ \forall^\infty k \ g(k) \neq f(k)\}$.}
\vskip 0.2truecm
\noindent
 From this, the Note and Theorem~1 we deduce
that ${\rm non}({\cal M})$ is the smallest cardinality of a family
$S\subseteq \omega^\omega$ with the property that for each
$f:\omega\rightarrow\bigcup_{n\in\omega}\omega^n$ there exists a
sequence
$\{a_n\}_{n\in\omega}$ belonging to $S$ such that for infinitely
many $i\in\omega$ the infinite sequence $\{a_{i+n}\}_{n\in\omega}$
extends the finite sequence $f(i)$.

Let ${\cal M}(\{0,1\}^\omega)$ denote the ideal of first category
subsets of the Cantor discontinuum $\{0,1\}^\omega$.
Obviously:
\begin{description}
\item[$(\ast\!\ast\!\ast)$]
\hspace{1.0cm} ${\rm non}({\cal M})
={\rm min}\{{\rm card}\ S:
S \subseteq \{0,1\}^\omega,\; S\not\in{\cal M}(\{0,1\}^\omega)\}$
\end{description}
\noindent
 From $(\ast\!\ast\!\ast)$, Theorem~1 and Theorem~2 we deduce that
${\rm non}({\cal M})$ is the smallest cardinality of a family
$S\subseteq\{0,1\}^\omega$ with the property that for each
$f:\omega\rightarrow\bigcup_{n\in\omega}\{0,1\}^n$ there exists a
sequence
$\{a_n\}_{n\in\omega}$ belonging to $S$ such that for infinitely
many $i\in\omega$ the infinite sequence $\{a_{i+n}\}_{n\in\omega}$
extends the finite sequence $f(i)$.
Another combinatorial
characterizations of ${\rm non}({\cal M})$ can be found in [5].
\par
\vskip 0.2truecm
{\bf Remark~3.} Errata to [5].
\begin{tabbing}
Page, line \= For \hspace{6.7cm}
\= Read\\
$22_{10}$ \> or all but finitely many \> for all but finitely
many \\
$26^{6}$ \> $k \in \omega (A \cup B)$ \> $k \in \omega \setminus
(A \cup B)$ \\
$27^{4}$ \> $C \Leftrightarrow C_1 \Leftarrow C_2 \Leftarrow
C_3 \Leftarrow ...$ \>\\
\hspace{3.0cm}
{$C \Leftrightarrow C_1 \Leftarrow C_2 \Leftarrow
C_3 \Leftarrow ...
\{f_{\{ a_n \} }: \{ a_n \}\in\Phi\}$ is unbounded}
\end{tabbing}

{\bf Acknowledgement.} The author is deeply indebted to Tomek
Bartoszy\'{n}ski for his communication [3].
\begin{center}
{\bf References}
\end{center}
\begin{enumerate}
\item T.~Bartoszy\'{n}ski, {\it Combinatorial aspects of measure
and category}, Fund. Math. 127 (1987), pp. 225-239.
\item T.~Bartoszy\'{n}ski and H.~Judah, {\it Set theory: on the
structure of the real line}, A.~K.~Peters Ltd., Wellesley MA 1995.
\item T.~Bartoszy\'{n}ski, {\it Private communication},
January 2000.
\item A.~W.~Miller, {\it A characterization of the least cardinal
for which the Baire category theorem fails}, Proc. Amer. Math.
Soc. 86 (1982), pp. 498-502.
\item A.~Tyszka, {\it On the minimal cardinality of a subset of
{\bf R} which is not of first category}, J. Nat. Geom. 17 (2000),
pp. 21-28.
\end{enumerate}
\begin{flushleft}
{\it Technical Faculty\\
Hugo Ko{\l}{\l}\c{a}taj University\\
Balicka 104, PL-30-149 Krak\'{o}w, Poland\\
rttyszka@cyf-kr.edu.pl\\
http://www.cyf-kr.edu.pl/\symbol{126}rttyszka}
\end{flushleft}
\end{document}